\magnification = \magstep1

\input amssym.def
\input amssym.tex
\def\iitem#1{\goodbreak\par\noindent{\bf #1}}
\def\Z{{\Bbb Z}}
\def\C{{\Bbb C}}
\def\Q{{\Bbb Q}}

\font\bbf=cmbx12

\def\qed{{~~~\vrule height .75em width .4em depth .2em}}
\def\irr#1{{\rm Irr}(#1)}
\def\ibr#1{{\rm IBr}(#1)}

\def\cent#1#2{{\bf C}_{#1}(#2)}

\def\syl#1#2{{\rm Syl}_#1(#2)}

\def\nor{\triangleleft\,}

\def\oh#1#2{{\bf O}_{#1}(#2)}

\def\zent#1{{\bf Z}(#1)}
\def\det#1{{\rm det}(#1)}
\def\ker#1{{\rm ker}(#1)}

\def\norm#1#2{{\bf N}_{#1}(#2)}

\let\phi=\varphi
\def\div{\raise 1pt \hbox{\big|}}
\def\ndiv{{\raise 1pt\hbox{/}\kern-4pt\div}}
\def\xdot{ \buildrel\textstyle\lower.5ex\hbox{.}\over\times}

\def\sbs{\subseteq}

\magnification = \magstep1

\font\bbf=cmbx12
\def\irr#1{{\rm Irr}(#1)}

\def\cent#1#2{{\bf C}_{#1}(#2)}

\def\syl#1#2{{\rm Syl}_#1(#2)}

\def\nor{\triangleleft\,}

\def\oh#1#2{{\bf O}_{#1}(#2)}

\def\zent#1{{\bf Z}(#1)}
\def\det#1{{\rm det}(#1)}
\def\ker#1{{\rm ker}(#1)}

\def\norm#1#2{{\bf N}_{#1}(#2)}

\let\phi=\varphi
\def\div{\raise 1pt \hbox{\big|}}
\def\ndiv{{\raise 1pt\hbox{/}\kern-4pt\div}}
\def\xdot{ \buildrel\textstyle\lower.5ex\hbox{.}\over\times}

{\nopagenumbers
\vglue 2 truein
\font\bbf = cmbx12
\centerline{{\bbf Irreducible characters taking root of unity values on $p$-singular 
elements}}
\bigskip
\centerline{by}
\medskip
\centerline{Gabriel Navarro}
\centerline{Departament d'\`Algebra}
\centerline{Universitat de Val\`encia}
\centerline{46100 Burjassot}
\centerline{SPAIN}
\smallskip
\centerline{E-mail: gabriel.navarro@uv.es}
\medskip
\centerline{and}
\medskip
\centerline{Geoffrey R. Robinson}
\centerline{Institute of Mathematics}
\centerline{University of Aberdeen}
\centerline{Aberdeen}
\centerline{AB24 3UE SCOTLAND}
\smallskip
\centerline{E-mail: g.r.robinson@abdn.ac.uk}
\vglue 1truein

ABSTRACT. In this paper we study finite $p$-solvable groups
having irreducible complex characters $\chi \in \irr G$ which
take roots of unity values on the $p$-singular elements of $G$.

\bigskip

\vfil

\eject}

\iitem{1. Introduction.}~~If $G$ is a finite group
and $p$ is a prime, in this paper we study
irreducible complex characters $\chi \in \irr G$
which take roots of unity values on the
$p$-singular elements of $G$.

There are several reasons to
study such characters.
In [NR], we gave a proof of a conjecture of J. Carlson, N. Mazza
and J. Th\'evenaz showing that an endo-trivial simple module of a finite $p$-solvable 
group 
$G$
with Sylow $p$-subgroups of rank at least 2 is one-dimensional. Recall that, as 
introduced
by E. C. Dade,
a simple $KG$-module $V$
over a field $K$ of characteristic $p$ is {\bf endo-trivial} if
$$V \otimes V^* \cong 1 \oplus P \, ,$$ where
$P$ is a projective module. (Later the second author
studied in [R2] $RG$-lattices satisfying this condition.) Inspired by this idea, we now 
turn our 
attention to the
set of complex irreducible characters $\irr G$ and consider 
characters $\chi \in \irr G$
such that
$$\chi \bar\chi=1 + \sum_i a_i\Phi_{\phi_i} \, ,$$
where $\Phi_{\phi_i}$ is the projective indecomposable character
of the irreducible Brauer character $\phi_i \in \ibr G$. Since $\Phi_{\phi_i}$ vanishes
on the $p$-singular elements of $G$, such a character $\chi$ satisfies that $|\chi(x)
|=1$
for every $p$-singular $x \in G$, and conversely.  (By elementary character theory
what we have in fact is that $\chi(x)$ is  a root of unity for every $p$-singular $x \in 
G$.
See Lemma (2.2) below.)

\medskip

Under certain
natural circumstances, the cyclic
defect theory provides  examples of these characters in groups
with cyclic Sylow $p$-subgroups.   It turns out that this condition also holds in 
several simple groups.
For instance, the irreducible
characters of $PSL(2,p^n)$ of degree $q-1$ take the value $-1$ on all its $p$-singular 
elements. Furthermore,
if $G$ is a simple group of Lie type of rank 1 in characteristic $r$, then $\cent Gx$ is 
an $r$-group for each non-trivial
$r$-element $x \in G$. Hence, if $\chi$ is the Steinberg character,
for each non-trivial $r'$-element $y$, we have $\chi(y) = \pm |\cent Gy|_r = \pm 1$.
In particular, for each prime $p$ other than $r$ which divides $|G|$,
we have $|\chi(y)|=1$ for each $p$-singular element $y$ of $G$.
(As pointed out to us by P. H. Tiep, if $q$ is a power of
the prime $p$, then $PSL_3(q)$,  $PSU_3(q)$,
$Sz(q)$  and $^2G_2(q)$, among others,  have also characters of this type.)

\medskip

Our interest in this paper is in $p$-solvable groups. In fact, Theorem A below
implies our module theoretic result from [NR].
\medskip

\iitem{THEOREM A.}~~{\sl Suppose that $G$ is a $p$-solvable group
with Sylow $p$-subgroups of rank at least $2$. Let $\chi \in \irr G$ be faithful
taking roots of unity values on the $p$-singular elements of $G$. If $\oh pG=1$,
then $\chi$ is linear.}

\medskip

If the Sylow $p$-subgroups of $G$ are cyclic or generalized quaternion, then there
are many examples showing that Theorem A is false (for instance, in certain
Frobenius actions). 

\medskip

Our result in Theorem A shows  that
under certain circumstances, if the values of an irreducible character on $p$-singular 
elements are
roots of unity, then the same is true everywhere; a fact that can be seen as another
example of character values on $p$-singular elements
controlling  representation theoretic invariants. 
\medskip

In the case where $\oh pG>1$, and without assuming $p$-solvability,
the situation is also quite tight.

\medskip

\iitem{THEOREM B.}~~{\sl   Suppose that $\chi \in \irr G$ is faithful, non-linear,
and such that it takes roots of unity values on the $p$-singular elements
of $G$. Assume that $P=\oh pG>1$. Then either $p$ is odd,
$P \in \syl pG$ is elementary abelian, $\cent GP=P \times \zent G$,  $G/\cent GP$
acts transitively on $\irr P -\{1_P\}$, and $\chi(1)=|P|-1$,
 or $p=2$ and $G=S_4 \times \zent G$ 
with 
$|\zent G|$ odd.}
\medskip
From Theorem B, we see that
$G/\zent G$ is  a doubly transitive
permutation group (in the action with $P$ acting by translation
on itself, and $G/\cent GP$ acting by conjugation on $P$).
Since the  finite 2-transitive groups are known, it is possible
to classify all the groups satisfying Theorem B,
but we have not
attempted this here.

\medskip

Finally, we mention the work of the second author in [R3], in which generalized 
characters
taking roots of unity values on all non-identity elements are studied.

\bigskip

\iitem{2. Proofs.}~~Let us start with the following.
\medskip

\iitem{(2.1) LEMMA.}~~{\sl Suppose that $\xi \in \C$ is a root of unity,
and let $e \in \Z$ be an integer. If $e$ divides $\xi$ in the ring of algebraic
integers, then $e=\pm 1$.}

\medskip

\iitem{Proof.}~~We have that $e\alpha=\xi$ for
some algebraic integer $\alpha$. Then
$1 =\xi \bar\xi=e^2\alpha\bar\alpha$,
and $(1/e^2)=\alpha\bar\alpha$ is a
rational algebraic integer. So $1/e^2 \in \Z$. \qed

\medskip

Our notation for characters, mainly follows [I].

\medskip

\iitem{(2.2) LEMMA.}~~{\sl Suppose that $\chi \in \irr G$. Then the following conditions 

are equivalent:
\smallskip

(i)~~The character  $\Psi=\chi\bar\chi-1_G$ vanishes on the
$p$-singular elements of $G$.

(ii)~~$|\chi(x)|=1$ for all $p$-singular $x \in G$.

(iii)~~$\chi(x)$ is a root of unity for all $p$-singular $x \in G$.}

\medskip

\iitem{Proof.}~~We know that $\Psi=\chi \bar\chi -1_G$ is
a character of $G$ since $1=[\chi,\chi]=[\chi\bar\chi,1_G]$. Now, if $\chi$ is any 
character
of $G$, $g \in G$ and $|\chi(g)|=1$ it is well-known that $\chi(g)$
is a root of unity. 
(This fact can be found as Problem (3.2) of [I], for instance, and we
sketch a proof for the interested reader.
Let $K$ be the cyclotomic field of $n$-th roots of unity,
of degree $\phi(n)$ over $\Q$. Let $\alpha \in K$ be an algebraic
integer with $|\alpha|=1$. Suppose that $\alpha$
has minimal polynomial
$p(x)= (x-\alpha_1)(x-\alpha_2) \cdots (x-\alpha_t) \in \Z[x]$, where the $\alpha_i$'s
are the Galois conjugates of $\alpha=\alpha_1$ by the  group ${\rm Gal}(K/\Q)$.
Hence, $t \le \phi(n)$. Now,
since ${\rm Gal}(K/\Q)$ is abelian, then Galois conjugation commutes
with complex conjugation, and therefore $|\alpha_i|=1$ for all $i$.
Now notice that the absolute value of the coefficient of $x^j$
in $p(x)$ is at most the binomial coefficient ${t \choose j}$,
so there only finitely many possibilities for $p(x)$. In particular,
we deduce that there is only a finite number
of algebraic integers $\alpha \in K$ with $|\alpha|=1$.
Therefore the set $\{ \alpha^m \, |\, m \ge 0 \}$ is finite,
and we conclude that $\alpha$ is a root of unity.)  \qed

\medskip

Although we shall not need it in our proofs,
we mention now that a character of a finite group
vanishes on the $p$-singular elements of $G$ if and only
if it is a $\Z$-linear combination of the projective
indecomposable characters $\Phi_{\phi_i}$ (by Corollary (2.17) of [N], for instance).

\medskip
We shall frequently use the following lemma.

\medskip

\iitem{(2.3) LEMMA.}~~{\sl  Suppose that $\chi \in \irr G$ takes roots
of unity values on every $p$-singular element $x \in G$. Then
\smallskip

(a)~~If $S \in \syl pG$, then $\chi(1)^2 \equiv 1$ mod $|S|$. Hence, if $p$
is odd, then
$\chi(1) \equiv \pm 1$ mod $|S|$. If $|S|=2^a$ and $a\ge 2$, then $\chi(1) \equiv \pm 1 
$ mod $2^a$
or $\chi(1) \equiv 2^{a-1} \pm 1$ mod $2^a$. In any case, $\chi(1)$ is not divisible by 
$p$.

(b)~~If $M \nor G$ has order divisible by $p$,
then $\chi_M$ is a sum of distinct irreducible characters of $M$.
In particular, if $\chi$ is non-linear and faithful, then
$\zent G$ is a $p'$-group.

(c)~~If $M$ is a subgroup of $G$ having a normal $p$-complement $N$
and a non-trivial Sylow $p$-subgroup $P$, then
$$[\chi_M, \chi_M]= 1 + {[\chi_N, \chi_N] -1 \over |P|} \, .$$}

\medskip

\iitem{Proof.}~~By Lemma (2.2),
we can write $\chi \bar\chi=1 + \Psi$, where $\Psi$ is a character of $G$ such that
$\Psi(g)=0$ whenever $g \in G$ is $p$-singular. If $S \in \syl pG$, then
we have that $\psi_S$ is a multiple of the regular character  $\rho_S$
of $S$,
so $\psi(1)$ is divisible by $|S|$. Then (a) follows.
To prove part
(b), using Clifford's theorem,
write
$\chi_M=e \sum_{i=1}^t \theta_i$, where the $\theta_i \in \irr M$ are distinct (and 
$G$-conjugate).
If $x \in M$ has order $p$, then
 $$\chi(x)=e(\sum_{i=1}^t \theta_i(x)) \, .$$
Hence, the root of unity $\chi(x)$ is divisible by $e$
in the ring of algebraic integers, and we apply Lemma (2.1).

Suppose finally that $M$ has a normal $p$-complement $N$ and
a Sylow $p$-subgroup $P>1$.  Then every $g \in M-N$ is $p$-singular.
Then
$$|M|[\chi_M,\chi_M]= \sum_{g \in M} |\chi(g)|^2=
\sum_{g \in M-N} |\chi(g)|^2 + \sum_{g \in N} |\chi(g)|^2=
(|M|-|N|) + |N|[\chi_N, \chi_N] \, ,$$
and part (c) easily follows.  \qed

\medskip
We shall also use the following lemma (due to the second author).
\medskip

\iitem{(2.4) LEMMA.}~~{\sl Suppose that $\theta$ is
a generalized character of an Abelian group $A$. If
$|\theta(a)|=1$ for all $a \in A-\{1\}$, then
$\theta=f\rho_A - e\mu$,
where $\rho_A$ is the regular character of $A$, $e$ is
a sign, and $\mu \in \irr A$.}

\medskip

\iitem{Proof.}~~This fact was first stated 
as a remark after Corollary 4 of [R1]. Another proof
is provided in the second paragraph of the proof of
Theorem (2.1) of [R3]. \qed

\medskip

Now, we dispose of Theorem B of the introduction.

\medskip

\iitem{(2.5) THEOREM.}~~{\sl   Suppose that $\chi \in \irr G$ is faithful, non-linear,
and such that it takes roots of unity values on the $p$-singular elements
of $G$. Assume that $P=\oh pG>1$. Then either
$P \in \syl pG$ is elementary abelian, $\cent GP=P \times \zent G$,  $G/\cent GP$
acts transitively on $\irr P-\{1_P\}$, and $\chi(1)=|P|-1$, or $p=2$ and
$G=S_4 \times \zent G$, where $|\zent G|$ is odd.}

\medskip

\iitem{Proof.}~~By Lemma (2.3.a), we know that $\chi$ has degree
not divisible by $p$. Hence, the
irreducible constituents of $\chi_P$
are linear. Since $\chi$ is faithful, we have that $P$ is abelian.
Now, by Lemma (2.3.b) and using that $\chi$
is faithful, then we can write
$\chi_P=\sum_{i=1}^t \lambda_i$, where $\lambda_i \in \irr P$ are
distinct and different from $1_P$.
Thus
$\chi(1) \le |P|-1$. Let $C=\cent GP$.
Since $[\chi_P, \lambda_1]=1$ and $P \sbs \zent C$,
we easily conclude that $\chi_C$ is a sum of
distinct linear characters of $C$. Using again that $\chi$ is faithful,
we conclude that $C$ is abelian.  Hence, we can write
$C=P \times N$, where $N=\oh {p'}G$. Now, write 
$$\chi_N=h\sum_{i=1}^r \mu_i \, ,$$
where $\mu_i \in \irr N$ are distinct.   Hence $\chi(1)=hr \le |P| -1$. Also,
$[\chi_C, 
\chi_C]=hr$.
Now, by Lemma (2.3.c), 
$$hr= 1 + {h^2 r -1\over |P|} \le 1 + {h^2 r -1\over hr+1}  \, .$$
From this inequality, we conclude that $r=1$ and that $\chi(1)=h=|P|-1$. Hence
$\chi_N=\chi(1)\mu_1$ and since $\chi$ is faithful and $\mu_1$ is linear,
we get $N\sbs \zent G$. Necessarily $N=\zent G$ (since $\zent G$ is a $p'$-group
by Lemma (2.3.b)). Also,
since $\chi(1)=|P|-1$, we have that $G/C$ acts faithfully
and  transitively on $\irr P-\{1_P\}$. In particular,
$P$ is an elementary abelian $p$-group.

So we may assume now that $P<S \in \syl pG$.
 If $p$ is odd, then by Lemma (2.3.a),
we have that $\chi(1) \ge   |S|-1$, where $S \in \syl pG$. Thus $|P|=|S|$ and  $P=S \in 
\syl pG$, which
is not possible. Hence, we have that $p=2$.
If $|P|=4$, then $G/C$ is necessarily $S_3$ and $G/\zent G=S_4$.
Since $|\zent G|$ is odd, by elementary group theory
we have that $G=\zent G \times S_4$.

Finally, assume that  $|P|>4$. By Lemma (2.3.a), we have that 
$|P|-1=\chi(1) \ge |S|/2 -1$ and  therefore $|S|=2|P|$. Since
$P$ is elementary abelian, it follows that the exponent of $S$
is at most 4.
Let $\lambda=\det \chi$, and let $\lambda_2$ be the 2-part of
$\lambda$. Now, let $\nu$ be a linear character of $G$ such that
$\nu^{\chi(1)}=\bar{\lambda_2}$, and notice that $\tau=\nu\chi$
is an irreducible character of $G$ with odd determinantal order
which takes roots of unity values on the $2$-singular elements
of $G$.  Since $\det \tau$ has odd order, it follows that
$\det \tau(x)=1$ on $2$-elements. Therefore
$\tau(x) \equiv \tau (1)$ mod 4 for every involution $x \in G$. 
Since $\tau(1)=|P|-1$ and $\tau(x)$ is a sign by hypothesis, we see that
$\tau(x)=-1$ for every involution $x \in G$.  Now let $U$ be
a cyclic subgroup of order 4. Since $\tau_U$
takes roots of unity values on the non-identity elements of $U$, it
follows by Lemma (2.4) that
$$\tau_U= f \rho_U -e\mu \, ,$$
where $e$ is a sign, $\mu \in \irr U$, and $\rho_U$ is the regular character of $U$.
Now, $|P|-1=4f -e$ and we deduce that $e=1$ and $f=|P|$.
Also, since $\det \tau_U=1$ and  $\tau_U= |P| \rho_U -\mu$,
we deduce that $\mu^{|P|-1}=1$ and hence $\mu=1$. It follows that
$\tau(x)=-1$ for every $x \in G$ of order 4.  Therefore
$\tau_S=a\rho_S -1_S$ and $\tau(1) \ge |S|-1$, which is the final contradiction. \qed

\medskip

We recall that a subgroup $H$ of a finite group $G$ is
said to be {\bf strongly $p$-embedded in $G$} if $p$ divides $|H|$ and $p$ does not
divide $|H \cap H^{g}|$ for each $g \in G -H$.

\medskip

The following is an essential part in the proof of Theorem A.
We use the Glauberman correspondence, and we refer the reader to [I], Chapter 13,
for a reference. 
But also, we shall use a well-known fact which follows from the Classification
of Finite Simple Groups, namely, that the outer automorphism 
group
of a  finite simple group of order prime to $p$ has cyclic Sylow $p$-subgroups (see, for 
example,
Theorem 7.1.2 of [GLS]).

\medskip

\iitem{(2.6) THEOREM.}~~{\sl Suppose that $G$ has a normal $p$-complement
$N$ and that $P \in \syl pG$ is abelian and contains
an elementary subgroup of order at least $p^2$. 
Let $\chi \in \irr G$ such that $\chi(x)$ is a root of unity for every $p$-singular $x 
\in G$.
Then $\chi(1)=1$.}

\medskip

\iitem{Proof.}~~We argue by induction
on $|G|$. By Lemma (2.3.a), we have that $\chi$ has $p'$-degree. Hence
$\chi_N=\mu \in 
\irr N$ (using Corollary (11.29) of [I]).
Now let $\hat\mu \in \irr G$ be the canonical extension of $\mu$ to $G$.
(This is the unique extension $\tau \in \irr N$ with determinantal order
$o(\tau)=o(\mu)$, using Corollary (6.28) of [I]). Then, by Corollary (6.17) of [I], we 
can write
$\chi=\hat\mu \lambda$ for some linear $\lambda \in \irr P$,
and therefore may assume that $\chi=\hat\mu$.

Since $(\hat\mu)_{N \times \oh pG}$ is by definition the
canonical extension of $\mu$ to $N \times \oh pG$,
we see that $(\hat\mu)_{N \times \oh pG}=\mu \times 1_{\oh pG}$.
Hence, if  $1\ne x \in \oh pG$, then  
$\hat\mu (x)=\mu(1)$ is a root of unity, and deduce that $\mu(1)=1$.
So we may assume that $\oh pG=1$.

Now we claim that $\hat\mu$ takes roots of unity
values on the $p$-singular
elements of $G$  if and only if for every $1 \ne x \in P$, the $\langle 
x\rangle$-Glauberman
correspondent of $\mu$ is linear. If $1 \ne x \in P$ and
$y \in \cent Nx$, then we know that
$$\hat\mu(yx)=\epsilon \mu^{*}(y) \, ,$$
where $\mu^{*}$ is the $\langle x\rangle$-Glauberman correspondent of 
$\mu$ and $\epsilon$ is a sign. (See Theorem (13.6) and (13.14) of [I].)
If $\mu^{*}$ is linear for every $1 \ne x \in P$, then we see that $\hat\mu$ takes
roots of unity values on $p$-singular elements.
Conversely, suppose that $\hat\mu$ takes roots
of unity values on $p$-singular elements.
If $1 \ne x \in P$ and $C=\cent Gx$, then 
$\hat\mu(x)=\epsilon\mu^{*}(1)$ is a root of unity, and then
 $\mu^{*}(1)=1$. This proves the claim.

Next we prove that $\chi$ is primitive. Suppose that $\tau^G=\chi$, where
$\tau \in \irr H$. Since $\chi$ has $p'$-degree, then $H$ contains
a Sylow $p$-subgroup of $G$, which we may assume is $P$. We claim that $H$ is
strongly $p$-embedded in $G$. Let $x \in H$
be a $p$-element. Since  $G$ has a normal $p$-complement
and abelian Sylow $p$-subgroups,
by elementary group theory we know that
 $x^g \in H$ for
some $g \in G$ if and only if $x^g=x^h$ for some $h \in H$.
Then 
$$\chi(x)=(1/|H|) \sum_{g \in G, x^g \in H} \tau(x^g)=
(1/|H|)m(x)\tau(x) \, ,$$
where $m(x)=|\{ g \in G \, |\, x^g \in H\}|$.
Now,   $x^g \in H$ if and only if $x^{gh} \in H$ for
every $h \in H$,  and we see that $m(x)/|H|$ is an integer dividing a root
of unity. Then $m(x)=|H|$ by Lemma (2.1). Hence if $x^g \in H$, then $g \in H$,
and therefore $H$ is strongly embedded in $G$, as claimed. In particular,
if $1 \ne u \in P$ and $z \in \cent Gu$, then $u \in H \cap H^z$ and therefore
$z \in H$.
Now, since $P$ contains as a
subgroup $C_p \times C_p$,
we may apply Theorem 6.2.4 of
[G], for example, to conclude that  $$N \sbs  \langle
\cent Gu \, |\, 1 \ne  u \in P \rangle  \sbs H\, $$ 
and we deduce that $H=G$. This proves our claim that $\chi$
is primitive.

Now, suppose that $U$ is a proper
normal subgroup of $N$ that admits $P$ and let $\theta \in \irr U$ be
$P$-invariant  under $\mu$. Let $1 \ne x \in P$, and let $\mu^* \in \irr{\cent Nx}$ be 
the
$\langle x \rangle$-Glauberman correspondent of $\mu$. By Theorem (13.29) of [I],
we have that $\theta^*$ lies
under $\mu^*$, and we conclude that $\theta^*$ is also
linear. Hence $\hat\theta \in \irr{UP}$
also takes roots of unity values on the $p$-singular
elements of $UP$, and by induction we have that $\theta$ is linear.

Now, let $N/Z$ be a chief factor of $G$, and let $\theta \in \irr Z$ be
$P$-invariant under $\mu$, which we know is linear.
Since $\hat\mu$ is primitive, then $\mu_Z$ is
a multiple of  $\theta$, so working
in $G/\ker \theta$, we may assume that $Z \sbs \zent G$
for every chief factor $N/Z$ of $G$.
Suppose that $N/Z$ is abelian, and let $1 \ne x \in P$.
Then $E/Z=\cent{N/Z}x \nor G/Z$ (using that $P$ is
abelian), and we deduce that $E=Z$
or $E=N$. If $E=N$, then $[N,x]=1$ by coprime action
and $x \in \oh pG=1$. Therefore $P$ acts Frobeniusly on $N/Z$,
but this is not possible because $P$ is not cyclic. So
we may assume that $N$ is perfect. 
Suppose that $N/Z$ is a simple group.
Since $\oh{p}{G/Z}=1$, it follows that $P$ is a subgroup
of ${\rm Aut}(N/Z)$, which we know by the Classification of Finite Simple Groups,
has cyclic Sylow $p$-subgroups. This is not possible.
Hence, 
$$N/Z= S_1/Z \times  \cdots \times S_k/Z$$ is a direct product
of non-abelian simple groups, where $k>1$, which are permuted by $P$,
where each $S_i$ is perfect and $[S_i,S_j]=1$.
Now, let $1 \ne x \in P$, and notice
that $\langle x \rangle$ transitively permutes the elements of the set
$\{S_1/Z, \ldots, S_k/Z \}$, because $\prod_{y \in \langle x \rangle} (S_i/Z)^y$
is normal in $G$.  Write $S_1=S$ and
define $\phi : S \rightarrow \cent Gx$ by
$$\phi(s)= s s^x \cdots s^{x^{p-1}} \, .$$ This is a well defined
homomorphism. If $s \in \ker \phi$,
then $s \in Z$ (using that $G/Z$ is a direct product) and $s^p=1$, so $s=1$.
Also, if $y \in \cent Gx$, then, using  again that $G/Z$ is a direct product,
we have that $y=ss^x \cdots s^{x^{p-1}}z$ for some $s \in S$ and $z\in Z$ .
Now $z=z_1^p$ for some $z_1 \in Z$ and $\phi(sz_1)=y$.
We conclude that $S \cong \cent Gx$, and in particular that $\cent Gx$ is
perfect. Since $\mu^* \in \irr{\cent Gx}$ is linear, then $\mu^*=1$ and $\mu=1$
by the uniqueness of the Glauberman correspondence.
Hence $\hat\mu$ is linear, as desired. \qed

\medskip

The following is Theorem A of the Introduction.

\medskip

\iitem{(2.7) THEOREM.}~~{\sl Let $G$ be a finite
$p$-solvable group. Suppose that $S \in \syl pG$ has
a non-cyclic subgroup of order $p^2$.
Suppose that $\chi \in \irr G$ is a faithful character that takes roots of unity values 
on the $p$-singular elements
of $G$. If  $\oh pG=1$,
then $\chi$ is linear.}

\medskip

\iitem{Proof.}~~Write $M=\oh {p'p}G=NP$, where $N=\oh {p'}G$ and $P$
is a $p$-group. Since $p$ divides $|G|$ and $G$ is $p$-solvable, we have that
$P>1$. Therefore $$\chi_{M}= \sum_{i=1}^s \chi_i$$
is a sum of different irreducible characters $\chi_i \in \irr M$, by Lemma (2.3.b).
Since $\chi$ has $p'$-degree, it follows that $(\chi_i)_N \in \irr N$ for all $i$.
Let us write
$$\chi_N=h \sum_{j=1}^r \mu_j \, ,$$
where $\mu_j \in \irr N$ are different. Notice that $s=hr$.
 By Lemma (2.3.c)
we obtain that
$$hr=s= [\chi_M, \chi_M]=1 + {h^2 r-1 \over |P|} \, .$$  
Hence $$(hr-1)|P|=(h^2 r -1)=h(hr-1) + (h-1)$$ and therefore
$hr-1$ divides $h-1$. Then $hr-1 \le h-1$, and we conclude that $r=1$.
So  $r=1$, $s=h$, and we have that
$\chi_N=h\mu$ for some $\mu \in \irr N$.
Since $(h-1)|P|=h^2-1$, then  
$h=1$ (and $\chi_M \in \irr M$) or $h=|P|-1$.

Let $\hat\mu \in \irr M$ be
the canonical extension of $\mu=\mu_1$ to $M$, so that $\hat\mu$ is
$G$-invariant. 

If $h=1$, then we have that $\chi_N \in \irr N$, and $\chi_M=\hat\mu \tau$
for some linear $\tau \in \irr P$. In this case, notice that
$\hat\mu (z)$ is a root of unity for every $p$-singular element
$z \in M$. 

  Suppose that $h=|P|-1$. We have that $\chi_i=\hat\mu \lambda_i$
for some unique $\lambda_i \in \irr{M/N}$, where the $\lambda_i$ are
distinct. By Clifford's theorem and using that $G=M\norm GP$,
we see that $\norm GP$ acts transitively on $\{ \lambda_1, \ldots, \lambda_s\}$.
Since $s=h=|P|-1$
then we conclude that $\irr{P} -\{1_P\}=\{ \lambda_1, \ldots, \lambda_s\}$.
In particular, $P$ is elementary abelian. If $P$ is cyclic, then
$G/M$ is abelian, and $P$ is a Sylow $p$-subgroup of $G$, which is not possible (
because
any Sylow $p$-subgroup of $G$ contains a subgroup of type $C_p \times C_p$).
Hence, $P$ contains at least a non-cyclic subgroup of order $p^2$. 
Also
$$\chi_M= \hat\mu \left(\sum_{\lambda \in \irr P -\{1_P\}} \lambda \right) \, .$$
Hence, if $z$ is a $p$-singular element of $M$, then
$\chi(z)=-\hat\mu(z)$
and we have that $\hat\mu$ takes roots of unity values
on $p$-singular elements.

In both cases, we apply Theorem (2.6) to conclude that
$\mu$ is linear. Since $\chi$ is faithful
and $\chi_N=\chi(1)\mu$, we conclude
that $\mu$ is faithful, so  $N=\oh{p'}G \sbs \zent G$. 
Now,  if $N<G$, then let $K/N$ be a chief factor of $G$.
Since $N=\oh {p'}G$, it follows that $K/N$ is
a $p$-group. Since $N \sbs \zent G$, then $K=N \times \oh pK$,
and this is impossible since $\oh pG=1$. Therefore,
we have that $G$ is a $p'$-group, and this is a final
contradiction. \qed

\medskip

We have mentioned in the Introduction that our present Theorem A
implies the module theoretic result in [NR]. In fact, if $V$ is a simple endo-trivial
$KG$-module of a finite $p$-solvable
group $G$ with $\oh pG=1$, $R$ is the usual complete local ring with $R/{\bf J}(R)=K$,
then by using the Fong-Swan theorem, it is possible
to lift $V$ to an $RG$-module $U$ that affords an irreducible $\chi \in \irr G$.
Now, by  Lemma 1 of [T], it easily follows that $U \otimes U^* \cong 1 \oplus P$,
where $P$ is a projective $RG$-module. Then $\chi(g)\bar\chi(g)-1=0$
for every $p$-singular $g \in G$ (using IV.2.5 in [F]), and we can apply Theorem A to
$\chi$.

\bigskip

\bigskip\bigskip

\centerline{REFERENCES}
\parindent = 0pt\frenchspacing
\medskip

\medskip

[F]~~W. Feit, {\sl The Representation Theory of Finite Groups}, North Holland, 1982.

\smallskip

[G]~~D. Gorenstein, {\sl Finite Groups}, Second edition, Chelsea Publishing Co., New
York, 1980.

\smallskip

[GLS]~~ D. Gorenstein, R. Lyons and R. Solomon, {\sl The Classification of the 
Finite Simple Groups, Number 3},
American Mathematical Society Mathematical Surveys and Monographs, Volume 40, 
Number 3, AMS, Providence, 1998.

\smallskip

[I]~~M. Isaacs, {\sl Character Theory of Finite Groups},   AMS Chelsea Publishing, 
Providence,
RI, 2006.

\smallskip

[N]~~G. Navarro, {\sl Characters and Blocks of Finite Groups}, Cambridge University 
Press, 1998.

\smallskip

[NR]~~G. Navarro, G. R. Robinson, On endo-trivial modules for $p$-solvable groups,
to appear in Math. Z. 
\smallskip

[R1]~~G. R. Robinson, A bound on norms of generalized characters with applications, J. 
Algebra {\bf 212} (1999), 660-668.

\smallskip

[R2]~~G. R. Robinson, Endotrivial irreducible lattices, submitted.

\smallskip
[R3]~~G. R. Robinson, Generalized characters whose values on non-identity elements are
roots of unity, submitted.

\smallskip

[T]~~J. G. Thompson, Vertices and sources, J. Algebra {\bf 6} (1967), 1-6.

\medskip

\bye